\documentclass[11pt]{article}
\usepackage{latexsym,amssymb,amsmath,amsthm,graphicx,float}

\newcommand{\R}{\mathbb{R}}

\newcommand{\<}{\langle}
\renewcommand{\>}{\rangle}

\newcommand{\eps}{\varepsilon}

\newtheorem{definition}{Definition}

\newtheorem{theorem}{Theorem}
\newtheorem{lemma}{Lemma}

\newtheorem{remark}[subsection]{Remark}

\title{Conditioning bounds for traveltime tomography in layered media}
\author{Hyoungsu Baek, Laurent Demanet \\ Department of Mathematics MIT \\ 77 Massachusetts Avenue Cambridge, MA  02139 USA}

\begin{document}
\maketitle

\begin{abstract}
This paper revisits the problem of recovering a smooth, isotropic, layered wave speed profile from surface traveltime information. While it is classic knowledge that the diving (refracted) rays classically determine the wave speed in a weakly well-posed fashion via the Abel transform, we show in this paper that traveltimes of reflected rays do not contain enough information to recover the medium in a well-posed manner, regardless of the discretization. The counterpart of the Abel transform in the case of reflected rays is a Fredholm kernel of the first kind which is shown to have singular values that decay at least root-exponentially. Kinematically equivalent media are characterized in terms of a sequence of matching moments. This severe conditioning issue comes on top of the well-known rearrangement ambiguity due to low velocity zones. Numerical experiments in an ideal scenario show that a waveform-based model inversion code fits data accurately while converging to the wrong wave speed profile.
\end{abstract}

{\bf Acknowledgments.} The authors would like to thank Guillaume Bal, Sergey Fomel, and William Symes for interesting discussions. This work was supported by a grant from Total SA. LD also acknowledges generous funding from the Alfred P. Sloan foundation and the National Science Foundation.

\section{Introduction}

\subsection{Problem setup and context}\label{sec:setup}

We consider the ray-theoretic limit of high-frequency waves propagating in a slab $0 \leq z \leq h$, made of a heterogeneous layered medium with smooth isotropic wave speed $c(z)$. We assume that waves can only be sent from, and recorded at the surface $z = 0$. Without loss of generality the waves are assumed to originate from the origin $x = z = 0$, as all points are equivalent on the surface. The transverse coordinate $x$ is assumed to be one-dimensional, as otherwise the problem would be radially symmetric about $x = 0$. We also assume that all other physical parameters that may affect wave dynamics, such as density, are constant.

\bigskip

The information available for the inversion is the traveltime $\tau$ of the various waves as a function of the recording position $x$. The two types of waves in a layered slab are
\begin{figure}[H]
\begin{minipage}{.3\linewidth}
\includegraphics[width=5.5cm]{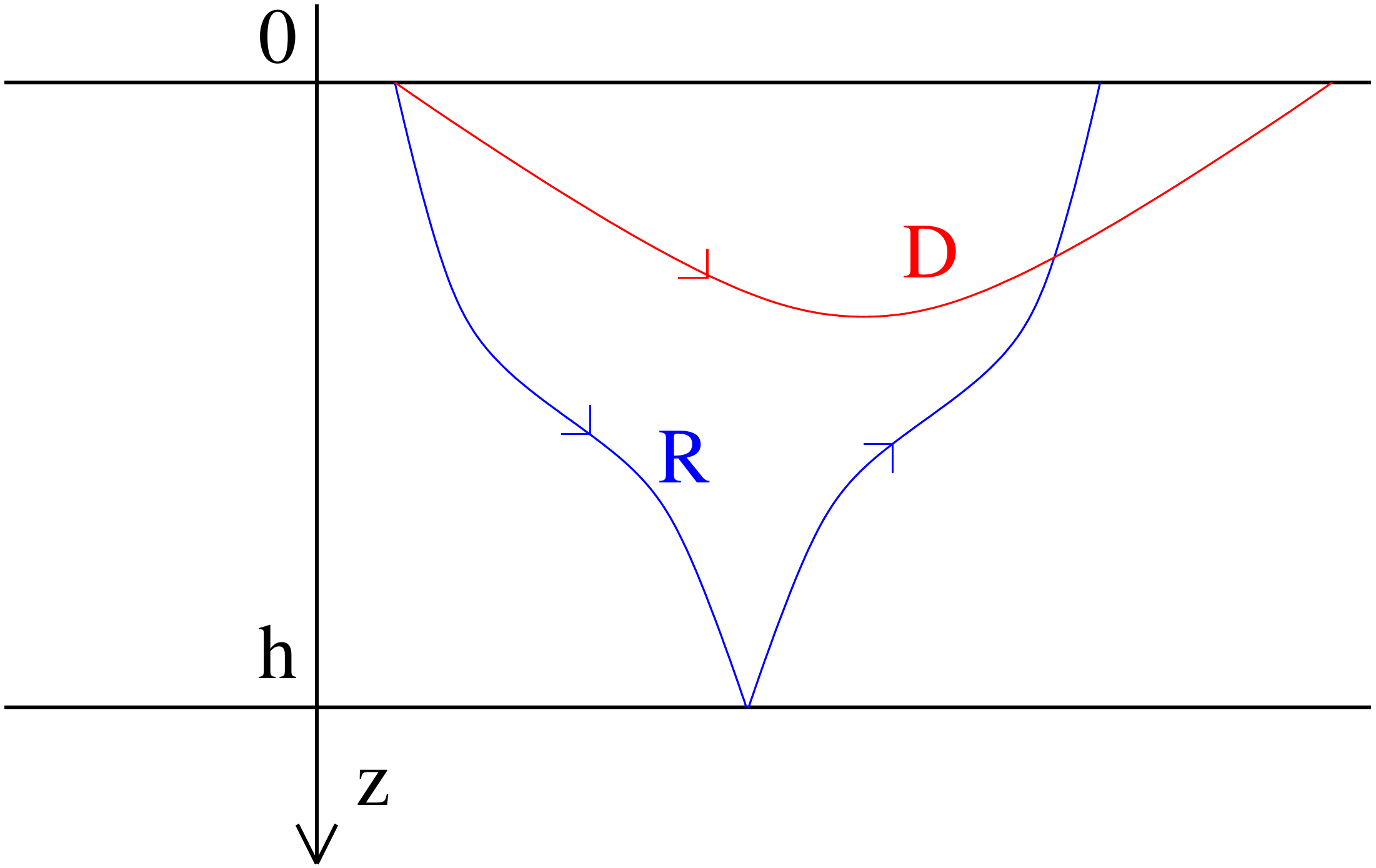}
\end{minipage}
\hfill
\begin{minipage}{.68\linewidth}
\begin{itemize}
\item \emph{diving}, or refracted waves, which arrive back at $z = 0$ from overturning before reaching $z = h$; and
\item transmitted waves, which arrive at $z = h$. The first \emph{reflected} wave, recorded at $z = 0$ after reflecting off of the boundary $z = h$, arrives twice later and twice farther than the transmitted wave, hence contains the same information. Multiply reflected waves also do not carry any new information.
\end{itemize}
\end{minipage}
\end{figure}

Diving waves occur for instance when $c(z)$ is monotonically increasing. We assume for simplicity that the type of a wave (diving or reflected) is a priori known. Waves that do not reflect (such as diving waves) are usually called ``transmitted" in the geophysics literature, so the word ``transmitted" is used very sparingly in the sequel to avoid confusion.

\bigskip

The inverse problem of recovering $c(z)$ from the traveltime information of \emph{diving} waves was solved circa 1910 by Herglotz \cite{Herglotz}, Wiechert and Geiger \cite{Wiechert}, and Bateman \cite{Bateman} in what is perhaps the first contribution by mathematicians to seismology. Their explicit formula takes the form of an inverse Abel transform and has been textbook material for a long time \cite{AkiRichards, Shearer, Bal,  Nowack1990, Nowack1997}. It will be reviewed in this paper, along with the analysis of its stability.

\bigskip

What can be said about the corresponding inverse problem for \emph{reflected waves}? 
Many authors have argued that this problem is quite different from that concerning diving waves. Firstly, there may not be an explicit formula to solve the problem. But more importantly, the problem has a completely different stability behavior. Qualitative discussion of ill-posedness of traveltime tomography was an active topic in the geophysics community in the late 1980s and early 1990s, see for instance Stork and Clayton \cite{Stork, ClaytonStork}; Bube, Langan, and Resnick \cite{Bube-1984, Bube-SIAM-1995, Bube-1995}; Ivansson \cite{Ivansson}; and Delprat-Jannaud and Lailly \cite{Lailly1, Lailly2}.  This paper aims to settle in a quantitative manner that, regardless of the discretization, there is not enough information in traveltime data from reflected rays to solve for a velocity profile $c(z)$ in a well-posed manner.

\bigskip

It is clear that in the one-dimensional case of a ray traveling from $z = 0$ to $z = h$ (with $h$ known), there is a fatal obstruction to solving the inverse problem. The only datum is the traveltime $\tau = \int_0^h 1/c(z) dz$, hence two smooth profiles with the same slowness integral will be indistinguishable. One would think that passing to a multidimensional situation may allow to recover a well-posed problem by triangulation --- the option to observe a fixed scene from different angles --- but that is not the case. Perhaps surprisingly, the presence of rays with different take-off positions and angles only marginally improves the determination of the velocity profile, at least in the layered case.

\bigskip

Determining a velocity profile form traveltime data is a nonlinear problem. As we will see, both in the diving and the reflected case, the forward model can be split into the composition of two operations:
\begin{itemize}
\item a \emph{nonlinear} operation of mapping the velocity profile to its decreasing rearrangement,  which is invertible when the function is monotone but not otherwise; followed by
\item a \emph{linear} integral operator acting on the inverse of this rearrangement. It is that integral operator which is invertible and relatively well-conditioned (of Volterra type) in the diving case, but always ill-conditioned (of Fredholm type) in the reflected case.
\end{itemize}
The possible lack of invertibility of the nonlinear step is well-understood: geophysicists refer to the lack of (increasing) monotonicity of the velocity profile in $z$ as the ``presence of low-velocity zones" \cite{AkiRichards, Shearer}. The characterization of the conditioning of the linear step seems to be less well understood and is the subject of this paper.

\bigskip

It should be mentioned that the meager results in this paper are far from shedding adequate light on the bigger problem of solvability and well-posedness of the general traveltime tomography problem, also called boundary rigidity problem. Much progress was obtained on this question recently; see \cite{StefanovUhlmann} and other upcoming publications by Uhlmann et al. For us, settling the layered special case serves to explain some disturbing numerical results that were observed in the scope of finite-frequency inversion of the background velocity in an idealized layered seismic setup. What we originally thought should have been a simple test case has now revealed itself to be a pathological example that cannot be solved. We hope it is useful to record this observation for the benefit of the community. We present a numerical example to this effect in the last section of this paper.

\bigskip

Although all our arguments assume a layered model $v(z)$, it is clear that they have an equivalent formulation in the radially symmetric case $v(r)$ via the so-called Earth flattening transformation. This is the original setting for the Herglotz-Wiechert formula.

\bigskip

It is a nice coincidence that some of the mathematics reviewed or used in the proofs originates from the first half of the 20th century, and should be credited to such first-rate analysts as Herglotz, Bateman, Hardy, Littlewood, and Szeg\H{o}.

\subsection{Kinematics}

In this section we review the solutions of the Hamiltonian system of geometrical optics in a layered medium. This classical material is covered in many places, including at least \cite{Whitham, Shearer, AkiRichards, Bal}. Let $\mathbf{x} = (x,z)$ for position and $\mathbf{p} = (p_x, p_z)$ for slowness; then
\begin{align*}
\dot{\mathbf{x}}(t) &= c(\mathbf{x}(t)) \frac{\mathbf{p}(t)}{|\mathbf{p}(t)|}, \\
\dot{\mathbf{p}}(t) &= - \nabla c(\mathbf{x}(t)) |\mathbf{p}(t)|.
\end{align*}
Since $c$ only depends on $z$, it follows that horizontal slowness is conserved and equals
\[
p_x \equiv p = \frac{\cos \theta_0}{c_0},
\]
where $\theta_0$ is the take-off angle that the ray leaving from the origin makes with the surface $z = 0$, and $c_0$ is the wave speed there. We now slightly abuse notations and write $c(z)$ for the wave speed. The rest of the system can be solved by writing 
\[
\frac{1}{c^2(z(t))} = |\mathbf{p}(t)|^2 = p^2 + p_z^2(t),
\]
isolating $p_z(t) = \sqrt{1/c^2(z(t)) - p^2}$, and using this expression in the equation for $z(t)$ to obtain
\[
\dot{z}(t) = c(z(t)) \sqrt{1 - p^2 c^2(z(t))}.
\]
Solving this ODE by separation of variables gives the expression of the traveltime $\tau$ as a function of $z$ and $p$:
\begin{equation}\label{eq:tau}
\tau(z,p) = \int_0^z \frac{1}{v(z,p)} \, dz, \qquad v(z,p) \equiv c(z) \sqrt{1 - p^2 c^2(z)}.
\end{equation}
The handy notation $v(z,p)$ refers to the vertical velocity. Returning to the equation for $x(t)$, we get
\[
\dot{x}(t) = \frac{dx}{dz} \dot{z}(t) = p \, c^2(z(t)),
\]
hence the horizontal position of the ray as a function of $z$ and $p$ is
\begin{equation}\label{eq:x}
x(z,p) = \int_0^z \frac{p \, c^2(z)}{v(z,p)} \, dz.
\end{equation}

\bigskip

The formulas (\ref{eq:tau}) and (\ref{eq:x}) can be used as is for transmitted rays, by letting $z = h$ in the upper bound of the integrals. Reflected rays obey the same expressions with a leading factor of 2. For short, we write $\tau(p)$ and $x(p)$ when $z = h$.

\bigskip

The traveltime and arrival position of a diving (refracted) ray, however, are obtained by following the ray until it reaches a turning point and then returns to the surface $z = 0$. A ray will turn if $v(z,p) = 0$, i.e., if it reaches the first $z = Z(p)$ where
\[
c(Z(p)) = 1/p.
\] 
Then, for diving rays,
\begin{equation}\label{eq:diving}
\tau(p) = 2 \int_0^{Z(p)} \frac{1}{v(z,p)} \, dz, \qquad x(p) = 2 \int_0^{Z(p)} \frac{p \, c^2(z)}{v(z,p)} \, dz.
\end{equation}

\bigskip

Data normally come in the form of one or more functions $T_i(x)$ of the transverse position $x$, but let us now explain how to introduce $p$ in this picture. 
Regardless of whether the ray is diving or reflected, the take-off angle $\theta_0$ is by symmetry the same as the angle that the ray makes with the surface $z = 0$ when recorded there. Hence $p = \cos \theta / c$ at the arrival point as well. It follows that $p$ is the rate of change of traveltime as a function of $x$:
\[
p_i = T_i'(x),
\]
where $i$ indexes the branch of the possibly multivalued traveltime. By inverting this relation we get the (unique) function $x(p)$. In turn, we get $\tau(p) = T_i(x(p))$. The step of forming $x(p)$ from its inverse function(s) may be numerically complicated, but does not in principle suffer from ill-conditioning. Hence in this paper we assume that $\tau(p)$ and $x(p)$ are given.

\bigskip

It is a unique feature of layered media that the knowledge of $T_i(x)$ implies that of the full scattering relation, i.e., of the take-off slowness vector in addition to the traveltime of each ray.

\subsection{Diving rays and the slowness distribution function}

The next (classical) step in solving the inverse problem is to change variables in (\ref{eq:diving}). As long as $c(z)$ is an increasing function of $z$, the relation $c(Z(q)) = 1/q$ introduced in the previous section defines the unique inverse function $Z(q)$. We then consider the Jacobian of $z$ with respect to $q^2$:
\[
F(q) =  \biggl\lvert \frac{dz}{dq^2} \biggl\rvert = \frac{ \lvert  {Z}'(q) \rvert}{2q}.
\]

\bigskip

For diving rays, the two relations in equation (\ref{eq:diving}) become
\begin{equation}\label{eq:diving2}
\tau(p) = 2 \int_p^{p_0} \frac{1}{\sqrt{q^2 - p^2}} \, (q^3 F(q)) \, dq, \qquad x(p) = 2 \int_p^{p_0} \frac{1}{\sqrt{q^2 - p^2}} \, (pqF(q)) \, dq.
\end{equation}
The upper bound is $p_0 = 1/c_0$, which is indeed greater than $p$ when $c(z)$ is increasing. We see that $\tau(p)$ and $x(p)$ in principle carry the same information; from the mathematical prospective it is sufficient to focus on $\tau(p)$ to determine $F(q)$.\footnote{It would be foolish in practice to ignore the position $x(p)$. If $F(q)$ is determined from $\tau(p)$ alone, the deviation of $x(p)$ from its integral expression above is at least an important indication of how well the assumption of layered medium is satisfied.}

\bigskip

It was the contribution of Herglotz \cite{Herglotz}, Wiechert and Geiger \cite{Wiechert}, and Bateman \cite{Bateman} to recognize that either of these Volterra integral relations can be reduced to the Abel transform, hence can be inverted in an explicit manner. In terms of $\tau(p)$, for instance,
\begin{equation}\label{eq:inverseAbel}
F(q) = - \frac{1}{\pi q^2} \int_q^\infty \frac{1}{\sqrt{p^2 - q^2}} \frac{d \tau}{dp}(p) dp.
\end{equation}

\bigskip

Putting numerical considerations aside, the solution of the inverse problem is now clear: 1) determine $F(q)$ from equation (\ref{eq:inverseAbel}), and 2) find $c^{-1}(z)$ as the inverse function of $Z(p) = \int_{p}^{p_0} 2qF(q)dq$. 

\bigskip

Although it is not the main topic of this paper, we detail for completeness the stability properties of the Abel inversion formula (\ref{eq:inverseAbel}) in Section \ref{sec:diving}. In a nutshell, the Abel transform is a half-integral, hence its inverse is a half-derivative. As a result, the inverse Abel transform is very mildly ill-posed. It is bounded between Lipschitz spaces with orders differing by $1/2$, and its singular values correspondingly decrease like $n^{-1/2}$.

\bigskip

The reasoning leading to an integral over slowness $q$ can be extended to non-monotonous $c(z)$ if we understand that $F(q)$ is now the \emph{slowness distribution function} (SDF),
\begin{equation}\label{eq:SDF}
F(q;p) = \int_0^{Z(p)} \delta \left( q - \frac{1}{c(z)} \right) \, dz. 
\end{equation}
This formula should convey the idea of a ``continuous histogram". The integrated or cumulative version of $F$ is sometimes called the ``layered cake" representation by analysts. If $F(\cdot; p)$ is too singular, the integral over $q$ should be understood in the sense of Stieltjes. Note that the dependence of $F$ on $p$ did not all of a sudden become crucial; the upper bound can be increased without consequence in the zone where $c^{-1}(z) < p$. What matters is that the integral avoids large $z > Z(p)$ for which $c^{-1}(z)$ increases back to $p$ and beyond; those values of $c^{-1}$ should not be counted in the SDF.

\bigskip

The expressions (\ref{eq:diving2}) still hold with $F(q;p)$ in place of $F(q)$ in more general situations when $c(z)$ may not be monotonically increasing. Since $\tau(p)$ and $x(p)$ only determine $F$, the main obstruction to solving for $c(z)$ is clear: \emph{any two profiles $c(z)$ which have the same SDF $F(q;p)$ will give rise to the same data $\tau(p), x(p)$.}  In other words, if $c(z)$ is the solution of the inverse problem, so will any (smooth) rearrangement within the interval $[0, Z(p)]$. This \emph{rearrangement ambiguity} does not pose a problem when $c(z)$ is increasing, but does in case $c(z)$ decreases before increasing back to $1/p$ at depth $Z(p)$. This explains the remark on ``low-velocity zones" in Section \ref{sec:setup}.

\bigskip

In addition to this rearrangement ambiguity, let us keep in mind that limited angular coverage is another major reason why practical traveltime tomography falls outside the scope of the Abel inversion formula. In fact, it is plausible that the techniques developed in the next sections would also help quantify the extent of the ill-posedness of the limited-data diving ray tomography problem.

\subsection{Reflected rays}

\begin{figure}[H]
\centering
\includegraphics[width=10cm]{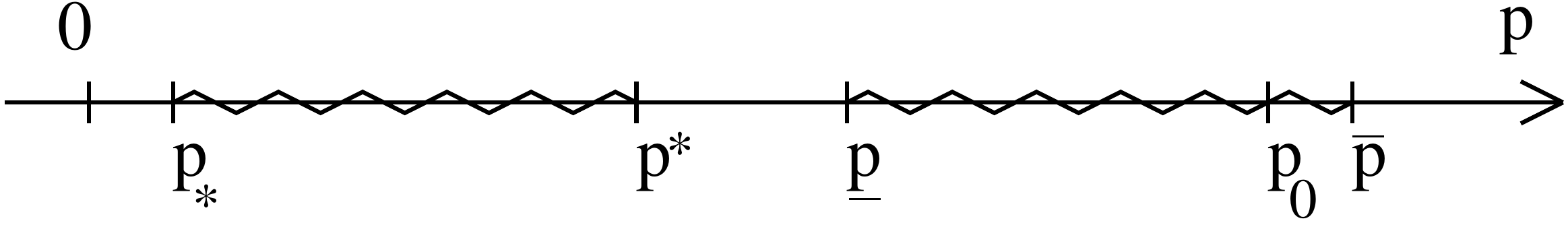}
\caption{Relative sizes of the remarkable slownesses, in the case of reflected rays. The interval $[\underline{p}, \overline{p}]$ is that of physical slownesses, i.e., $1/c(z)$ for some $z$. In particular, $p_0 = 1/c(0)$ belongs in this interval. The interval $[p_*, p^*]$ is that of observed horizontal slownesses, i.e., quantities of the form $p_0 \cos \theta$ where $\theta$ is the angle that the incident ray makes with the surface $z = 0$. Note that in the case of diving rays, the two intervals would overlap and $p^* > \underline{p}$.}\label{fig:p}
\end{figure}

For reflected rays, the situation is very different since we are back to formulas (\ref{eq:tau}) and (\ref{eq:x}) with $h$ in place of $z$ in the upper bound of the integrals. As a result, the $q$-integrals for reflected rays, in terms of the SDF introduced in the previous section, take the form
\begin{equation}\label{eq:reflected}
\tau(p) = 2 \int_{\underline{p}}^{\overline{p}} \frac{1}{\sqrt{q^2 - p^2}} \, (q^3 F(q;p)) \, dq, \qquad x(p) = 2 \int_{\underline{p}}^{\overline{p}} \frac{1}{\sqrt{q^2 - p^2}} \, (pqF(q;p)) \, dq.
\end{equation}
The integral bounds depend only on the medium properties; they are 
\[
\underline{p} = \min_{z \in [0,h]} \frac{1}{c(z)}, \qquad \overline{p} = \max_{z \in [0,h]} \frac{1}{c(z)}.
\]
The bounds $p_{*}$ and $p^{*}$ on the horizontal slowness variable $p$, on the other hand, relate to the availability of data. They are the slownesses for which the traveltime information of the reflected rays is known in the interval $[x(p_*), x(p^*)]$. The relations between the various remarkable slownesses are summarized in Figure \ref{fig:p}.

\bigskip

\bigskip

\emph{The lower bound $\underline{p}$ is now fixed, so we are dealing with first-kind Fredholm equations instead of Volterra equations for $F(q;p)$. This results in severe ill-conditioning of the inverse map in the case of reflected rays.}

\bigskip

The conditioning of the linear map in (\ref{eq:reflected}) depends on how $p^*$ relates to $\underline{p}$. This information is captured by the smallest angle $\psi$ that reflected rays make with the lines $z$ = const., whose cosine is
\[
\cos \psi = \frac{p^*}{\underline{p}}.
\] 
The larger this angle, the worse the conditioning of (\ref{eq:reflected}), i.e., the more unstable the inverse maps. But the problem is already ill-conditioned even if $\psi = 0$.

\bigskip

We will need the following functions of $\underline{p}, \overline{p}, p_*,$ and $p^*$: the contrast $e = (\overline{p} / \underline{p})^2$, and
\begin{align*}
\rho_* &= 1 + 2 \frac{1-(p_*/\underline{p})^2}{e-1} - \left[ \left( 1 + 2 \frac{1-(p_*/\underline{p})^2}{e-1} \right)^2 - 1 \right]^{1/2}, \\
\rho^* &= 1 + 2 \frac{1-(p^*/\underline{p})^2}{e-1} - \left[ \left( 1 + 2 \frac{1-(p^*/\underline{p})^2}{e-1} \right)^2 - 1 \right]^{1/2},
\end{align*}
and
\[
\alpha = \frac{1}{\rho^*}, \qquad \beta = \frac{\rho^* + \rho_* + 2}{\rho^* - \rho_*} + \left[ \left( \frac{\rho^* + \rho_* + 2}{\rho^* - \rho_*} \right)^2 - 1 \right]^{1/2}.
\]
It turns out that both $\alpha$ and $\beta$ are increasing functions of $\psi$ (when $p_*$ is fixed). Notice that $\alpha=1$, or $\rho^* = 1$, if and only if $\psi = 0$. We will see below that $\alpha$ gives rise to a lower bound $\sim \alpha^N$ on the condition number, while $\beta$ yields an upper bound $\sim \beta^N$, with $N$ yet to be defined. These formulas may look complicated, but there are numerical indications that the upper bound generated by $\beta$ is tight.

\bigskip

In the sequel we limit ourselves without loss of generality to the equation for $\tau(p)$. We let $A$ for the linear map from $F$ to $\tau$ in (\ref{eq:reflected}).

\bigskip

The notion of condition number is meaningful only if the number of degrees of freedom is limited. Most sampling schemes of interest would discretize the operator $A$ as a $M$-by-$N$ matrix. Specifically, let
\begin{itemize}
\item $\mathcal{P}_M$ be an orthogonal projector of rank $M$ on $L^2(p_*, p^*)$; and
\item $\mathcal{Q}_N$ be an orthogonal projector of rank $N$ on $L^2(\underline{p}, \overline{p})$. 
\end{itemize}
It is natural to consider $A_{MN} = \mathcal{P}_M A \mathcal{Q}_N$, and its condition number
\[
\kappa(A_{MN}) = \| A_{MN} \|_2 \; \| A^{+}_{MN} \|_2 = \frac{\sigma_{\max}(A_{MN})}{\sigma_{\min}(A_{MN})},
\]
where $+$ denotes the pseudo-inverse, and $\sigma$ denote the singular values.

\begin{definition}\label{def:reasonable}
A discretization is any couple $(\mathcal{P}_M, \mathcal{Q}_N)$ of orthogonal projectors. It is called \emph{reasonable for $A$} if
\[
\| \mathcal{P}_{M} A \mathcal{Q}_N \|_2 \geq \frac{1}{2} \| A \|_2.
\]
\end{definition}

\bigskip

\begin{theorem}\label{teo:main} Assume $p_* > 0$.  Let $\kappa_N$ be the condition number of $\mathcal{P}_M A \mathcal{Q}_N$.

\begin{itemize}
\item[(i)] For all discretizations $(\mathcal{P}_M, \mathcal{Q}_N)$, reasonable for $A$, and such that $\mathcal{Q}_N$ has rank $N$,
\begin{equation}\label{eq:lower1}
\kappa_N \geq C_\alpha N  \alpha^{N}.
\end{equation}
It also holds that (useful when $\alpha = 1$)
\begin{equation}\label{eq:lower2}
\kappa_N \geq C N^{-1/4} \, e^{(1.2465 \ldots) \sqrt{N}}.
\end{equation}
\item[(ii)] There exists a discretization $(\mathcal{P}_N, \mathcal{Q}_N)$, where both $\mathcal{P}_N$ and $\mathcal{Q}_{N}$ have rank $N$, such that
\begin{equation}\label{eq:upper}
\kappa_N \leq C_\beta \, N^{1/4} \, \beta^N.
\end{equation}
Furthermore, $\mathcal{Q}_N$ can be taken as the orthogonal projector on polynomials of $p^2$, of degree $N-1$.
\end{itemize}

\end{theorem}

\bigskip

We do not believe that the prefactors $N$, $N^{-1/4}$ and $N^{1/4}$ are sharp in any way. We give numerical indications in Section \ref{sec:num} that the rate $\beta^N$ of (\ref{eq:upper}) may be sharp, while the constants in (\ref{eq:lower1}) and (\ref{eq:lower2}) are not. That (\ref{eq:lower2}) should scale root-exponentially when $\alpha=1$ seems adequate, however.

\bigskip

It should be mentioned that the positivity constraint on $F(p)$ or $F(p;q)$ may be an important piece of information for solving the inverse problem in practice, but that it is does not call the linear stability analysis into question. As soon as $F(p)$ is strictly bounded away from zero on its support, no small perturbation can compromise positivity, hence the linearized perturbative theory applies.

\bigskip

Note that once $F(q;p)$ is found, there may still not be a unique $c(z)$ that corresponds to it.  The ambiguity of smooth rearrangements in the ``low velocity zones" is as much an issue for reflected rays as it is for diving rays. While the source of ill-conditioning is now twofold, we present a numerical example in Section \ref{sec:num} where the ill-conditioning detailed in Theorem \ref{teo:main} is in fact more problematic than the strict lack of uniqueness arising from the rearrangement ambiguity.

\subsection{Small $p$ asymptotics}\label{sec:asymp-intro}

We know from the previous section that there is a wide range of kinematically near-equivalent velocity profiles if one only considers data from reflected rays. It is possible to describe this range of velocities quite well in the case of small $p$, i.e., small offset between source and receiver. 

\bigskip

When $p=0$ (rays perpendicular to the layering), the problem is one-dimensional and the data reduce to the single number $\tau = \int_0^h c^{-1}(z) \, dz$. Any two smooth velocity profiles $c_1$, $c_2$ such that
\[
\int_0^h c_1^{-1}(z) \, dz = \int_0^h c_2^{-1}(z) \, dz
\]
will appear indistinguishable.

\bigskip

The extension of this observation to the case of small $p \geq 0$ is that traveltimes will be nearly equal provided the integrals of the \emph{odd powers} of $c$ are identical, namely
\[
\int_0^h c_1^{2n-1}(z) \, dz = \int_0^h c_2^{2n-1}(z) \, dz, \qquad \mbox{for } 0 \leq n < d
\]
for some (small) integer $d>0$. Then the traveltimes for $p \in [0, p^*]$ will match up to a remarkably small $O({(p^*)}^{2d})$. A justification is given in Section \ref{sec:asymp}.  The linearized version of the conditions above was found by Ivansson \cite{Ivansson}. See also the paper \cite{Bube-SIAM-1995} by Bube for a more extensive study of the slowness nullspace in the linearized regime.

\section{Theory for diving rays}\label{sec:diving}

In this section, we briefly review the Herglotz inversion formula (\ref{eq:inverseAbel}) and what is known about its stability properties.

\bigskip

It is convenient to change variables as $x = p^2$, $y = q^2$ and express (\ref{eq:diving2}) via an operator $A$ as
\begin{equation}\label{eq:newA-diving}
g(x) = A f(x) = \int_{x}^{\infty} \frac{1}{\sqrt{y-x}} f(y) \, dy, \qquad x \geq 0,
\end{equation}
with $g(x) = \tau(\sqrt{x})$ and $f(y) = y \, F(\sqrt{y})$. This relationship between $f$ and $g$ is (a version of) the Abel transform, or Abel equation \cite{Polyanin}. It is also called a Riemann-Liouville integral.

\bigskip

The key to the inversion formula for $A$ is to notice that it is an operator of fractional differentiation of order $1/2$ on the half-line. Indeed,
\[
A^2 f(x) = \int_x^\infty k(x,y) f(y) \, dy,
\]
with
\begin{align*}
k(x,y) &= \int_x^y \frac{1}{\sqrt{(z-x)(y-z)}} \, dz, \\
&= - \arctan \left( \frac{x+y-2z}{2 \sqrt{(z-x)(y-z)}} \right) \Bigg\vert_x^y, \\
&= \pi.
\end{align*}
As a result, $A^2 \, \frac{d}{dx} = - \pi$, hence
\[
A^{-1} = - \frac{1}{\pi} A \, \frac{d}{dx}.
\]
The inversion formula (\ref{eq:inverseAbel}) follows.

\bigskip

One can also quickly notice that 
\[
A^2 \, e^{sx} = - \frac{\pi}{s} e^{sx}, \qquad \mbox{Re}(s) < 0,
\]
from which one can (correctly) infer that
\begin{equation}\label{eq:AbelLaplace}
A \, e^{sx} = \left( - \frac{\pi}{s} \right)^{1/2} e^{sx}, \qquad \mbox{Re}(s) < 0.
\end{equation}
This leads to the well-known fact that $A$ is diagonal in the Laplace domain, that its powers form a semi-group, and that $A^{-1}$ can also be computed via the scaling $\sqrt{-s/\pi}$ in the Laplace domain \cite{Polyanin}. This procedure is not advisable numerically due to the ill-conditioning of the inverse Laplace transform.

\bigskip

Equation (\ref{eq:AbelLaplace}) also carries the information that any $z$ such that Re$(z) > 0$ is an eigenvalue of $A$, with square-integrable eigenfunction. The spectral theory of nonnormal operators such as $A$ is however quite pathological, so this observation is rather useless.  A natural finite dimensional approximation of $A^2$ would be a highly defective upper-triangular matrix with constant entries on and above the diagonal. Eigenvalues are not the right tool to study stability under inversion for such nonnormal operators or matrices.

\bigskip

Singular values, however, are perfectly informative for stability. The following result gives an explicit singular value expansion of $A$ in the illustrative case when $A$ acts on functions supported in $x \in [-1,1]$ -- otherwise some rescaling needs to be done. It should be credited to Johnstone and Silverman who proved a very similar result in \cite{Johnstone}.

\begin{theorem} (Johnstone-Silverman)
Assume $f \in C([-1,1])$. Then
\begin{equation}\label{eq:SVD-A}
A f (x) = \sum_{n} u_n(x) \sigma_n \< v_n, f \>_{\scriptsize r},
\end{equation}
with
\begin{align*}
u_n(x) &= (1-x^2)^{1/2} \, U_n(x), \qquad \mbox{(} U_n \mbox{ are Chebyshev polynomials of the second kind,)} \\
v_n(x) &= \left( \frac{n+1}{2} \right)^{1/2} \, P_n^{(0,1)}(x), \qquad \mbox{(} P_n^{(0,1)} \mbox{ are Jacobi polynomials,)} \\
\sigma_n &= \left( \frac{n+1}{2} \right)^{-1/2},
\end{align*}
and the inner product is
\[
\< f, g \>_{\scriptsize r} = \int_{-1}^1 f(x) g(x) (1+x) \, dx.
\]
\end{theorem}

\begin{remark}
Notice that the $v_n$ are an orthonormal basis for the ``right"  inner product $\< \cdot, \cdot \>_{r}$, whereas $u_n$ are an orthonormal basis for the ``left" inner product
\[
\< f, g \>_{\scriptsize \ell} = \int_{-1}^1 f(x) g(x) \frac{2}{\pi} \left( \frac{1+x}{1-x} \right)^{1/2} \, dx.
\]
The particular values of the $\sigma_n$ depend on the normalization of the inner products, but their decay rate $\sim n^{-1/2}$ does not.
\end{remark}

\begin{proof}
One first establishes that $A v_n = \sigma_n u_n$, from which (\ref{eq:SVD-A}) follows by orthonormality and completeness of the $v_n$.  The proof is a matter-of-fact induction argument which combines equation 22.13.11 in \cite{Abramowitz} with the relations $(2n+1) P_n^{(0,1)} = (n+1) P_n + n P_{n-1}$; $ \, T_n = U_n - x U_{n-1}$; and $U_{n+1} = 2x U_n - U_{n-1}$. Johnstone and Silverman claim that there is a less artificial way of obtaining relations such as $A v_n = \sigma_n u_n$ via hypergeometric functions.
\end{proof}

\bigskip

Since the singular values $\sigma_n$ decay like $n^{-1/2}$, so will the singular values of any good discretization of $A$. As a result, we can expect that a matrix discretizing $A$ on $N$ points would have a $O(\sqrt{N})$ condition number. This qualifies as very mild ill-posedness.

\bigskip

We may also understand the stability properties of $A^{-1}$ through boundedness estimates in adequate functional spaces. Consider the Lipschitz space $\mbox{Lip}(\alpha)$ of functions with bounded $\alpha$ semi-norm \cite{DeVore}
\[
\| f \|_\alpha = \left\{ \begin{array}{ll}
         \sup_{x \ne y} \frac{|f(x) - f(y)|}{|x-y|^\alpha} & \mbox{if $0 < \alpha \leq 1$};\\
         \| f^{(\lfloor \alpha \rfloor)} \|_{\alpha - \lfloor \alpha \rfloor}& \mbox{if $\alpha > 1$}.\end{array} \right. 
\] 

Hardy and Littlewood studied boundedness of fractional integration on Lipschitz spaces \cite{HardyLittlewood}. Their conclusion for $A$, properly modernized, is that for all $f \in \mbox{Lip}(\alpha)$, $\alpha > 0$,
\[
\| A f \|_{\alpha+1/2} \leq C \| f \|_{\alpha}, \qquad \| A^{-1} f \|_{\alpha - 1/2} \leq  C \| f \|_{\alpha}.
\]
This result again showcases the mild ill-posedness of inverting $A$.

\section{Theory for reflected rays}\label{sec:reflected}

In this section we prove Theorem \ref{teo:main}. As in the previous section we change variables as $x = p^2$, $y = q^2$, to obtain
\begin{equation}\label{eq:newA}
f(x) = \int_{\underline{x}}^{\overline{x}} \frac{1}{\sqrt{y-x}} g(y) \, dy, \qquad x \in [x_*, x^*],
\end{equation}
with
\[
f(x) = \tau(\sqrt{x}), \qquad g(y) = y \, F(\sqrt{y}; \sqrt{x}), \qquad [ \underline{x}, \overline{x}] = [\underline{p}^2, \overline{p}^2], \qquad [ x_*, x^*] = [p_*^2, {p^*}^2].
\]
The relationship between the bounds is the same as earlier, namely $x_* < x^* \leq \underline{x} < \overline{x}$. Since all $x$ and $y$ are bounded away from zero (from the assumption $p_* > 0$), considering the linear map between $f$ and $g$ --- rather than that between $\tau$ and $F$ --- changes the condition number by a factor independent of $M$, $N$. Hence it suffices to prove the claims of the theorem for (\ref{eq:newA}). We overload notations and reuse the letter $A$ for $f = A g$ as defined by (\ref{eq:newA}). Note that this equation is quite different from (\ref{eq:newA-diving}).

\bigskip

The singular values of $A$ are the square roots of the eigenvalues of $A^* A$,
\[
A^* A f(y') = \int_{\underline{x}}^{\overline{x}} k(y',y) f(y) \, dy, \qquad k(y',y) = \int_{x_*}^{x^*} \frac{1}{\sqrt{(y'-x)(y-x)}} \, dx.
\]
The kernel integrates to $k(y',y) = - 2 \log (2(\sqrt{y'-x} + \sqrt{y-x})) \vert_{x_*}^{x^*}$, which is clearly Hilbert-Schmidt on $[\underline{x},\overline{x}]^2$ even in the case when $\underline{x} = x^*$. Hence $A^* A$ is a compact operator. As a consequence of the general theory, it has a discrete set of eigenvalues (with square-integrable eigenvectors) which can only accumulate at the origin.

\subsection{Legendre expansion of the kernel}

A key to understanding the spectrum of $A^* A$ is that $1/\sqrt{y-x}$ has an explicit expansion in terms of the Legendre polynomials rescaled to the interval $[\underline{x}, \overline{x}]$. Consider the new variables
\[
\Sigma = \frac{\overline{x}+\underline{x}}{2}, \qquad \Delta = \frac{\overline{x}-\underline{x}}{2}, \qquad \tilde{x} = \frac{\Sigma - x}{\Delta}, \qquad \tilde{y} = \frac{\Sigma - y}{\Delta}.
\]
If $P_n(\tilde{y})$ denotes Legendre polynomial of degree $n$ with $\tilde{y} \in [-1,1]$, then
\[
p_n(y) = \sqrt{\frac{n+1/2}{\Delta}} P_n \left( \tilde{y} \right), \qquad n \geq 0
\]
is an orthonormal basis for $[\underline{x}, \overline{x}]$ with measure $dy$. 

\bigskip

The desired expansion is
\begin{equation}\label{eq:Legendre}
\frac{1}{\sqrt{y-x}} = \sqrt{\frac{2}{\Delta}} \, \sum_{n \geq 0} \rho^{n+1/2} P_n(\tilde{y}),
\end{equation}
where $\rho \leq 1$ is related to $\tilde{x} \geq 1$ through
\begin{equation}\label{eq:elliptical}
\tilde{x} = \frac{\rho + \rho^{-1}}{2} \qquad \Leftrightarrow \qquad \rho = \tilde{x} - \sqrt{\tilde{x}^2-1}.
\end{equation}
Equation (\ref{eq:Legendre}) is a straightforward consequence of the fact that  $1/\sqrt{1 - 2\tilde{y} \rho + \rho^2}$ is the generating function of the Legendre polynomials $P_n(\tilde{y})$. Equation (\ref{eq:elliptical}) is part a change of variables to and from elliptical coordinates in the complex plane; $\rho^{-1}$ is the elliptical radius of the Bernstein ellipse passing through $\tilde{x} \geq 1$, with foci at $\pm 1$. Ultimately, it is well-known that the speed of convergence of a series like (\ref{eq:Legendre}), or of the corresponding Chebyshev series, depends on the distance of the singularity at $\tilde{x}$ to the interval $[-1,1]$ in the complex plane.

\bigskip

The rescaled polynomials $p_n(y)$ provide a unitary change of basis for $A^* A$. Using equation (\ref{eq:Legendre}), it suffices to find the eigenvalues of the semi-infinite matrix
\begin{align*}
K_{m,n} &\equiv  \int_{\underline{x}}^{\overline{x}} \!  \int_{\underline{x}}^{\overline{x}} p_m(y') k(y',y) p_n(y) \, dy' dy', \\
&= 2 \, \Delta \left( (n+\frac{1}{2})(m+\frac{1}{2}) \right)^{-1/2} \int_{\tilde{x}_*}^{\tilde{x}^*} \left[ \, \tilde{x} - \sqrt{\tilde{x}^2 - 1} \,\, \right]^{m+n+1} \, d\tilde{x}, \qquad m, n \geq 0.
\end{align*}
Further passing to the $\rho$ variable via (\ref{eq:elliptical}), it follows that $K_{m,n}$ is (up to the normalization factor) a Hankel matrix of moments:
\begin{equation}\label{eq:Hankel}
K_{m,n} = \Delta \left( (n+\frac{1}{2})(m+\frac{1}{2}) \right)^{-1/2} \int \rho^{m+n} d\mu(\rho),
\end{equation}
with density
\[
\mu'(\rho) = \rho^{-1} - \rho, \qquad \rho \in [\rho_*, \rho^*],
\]
where the bounds $\rho_*$ and $\rho^*$ relate to $\tilde{x}^*$ and $\tilde{x}_*$ respectively through (\ref{eq:elliptical}); in particular $\rho^* = \tilde{x}_* - \sqrt{\tilde{x}_*^2 - 1}$. Note that $\rho_* < \rho^* \leq 1$.

\bigskip

We now address the decay of the eigenvalues of $K$ and of its finite-dimensional sections.

\subsection{Coarse lower bound (\ref{eq:lower1}) on the condition number}\label{sec:coarse-lower}

In this section we start by assuming $\rho^* < 1 \iff \psi > 0$, i.e., the rays never become horizontal and the kernel $k(y,y')$ is bounded.

\bigskip

The following two elementary lemmas detail how to deal with finite-dimensional projections of compact operators. Their proofs are nice homework exercises involving the Courant-Fischer min-max principle. In what follows eigenvalues are sorted in decreasing order, and projectors are considered on the domains over which they make sense.

\begin{lemma}\label{lem:one}
Let $\mathcal{P}_M$ and $\mathcal{Q}_N$ be two orthogonal projectors. Then for all $j \geq 1$,
\[
\lambda_j(\mathcal{Q}_N A^* \mathcal{P}_M A \mathcal{Q}_N) \leq \lambda_j(A^* A).
\]
\end{lemma}

\begin{lemma}\label{lem:two}
Let $\mathcal{R}_N$ be an orthogonal projector of rank $N$. Then
\[
\lambda_N(A^*A) \leq \lambda_1((I-\mathcal{R}_N) A^* A (I-\mathcal{R}_N)).
\]
\end{lemma}

Given our $N$-by-$N$ matrix $A_{M,N} =  \mathcal{P}_M A \mathcal{Q}_N$, where $ \mathcal{P}_M$ and $\mathcal{Q}_N$ are arbitrary, we use the two lemmas above to obtain the bound
\[
\lambda_N(A^*_{M,N} A_{M,N}) \leq \lambda_1((I-\mathcal{R}_N) A^* A (I-\mathcal{R}_N)).
\]
We choose $\mathcal{R}_N$ to be the orthogonal projector onto polynomials of degree $N \! - \! 1$ in $[-1,1]$, i.e. $(I-\mathcal{R}_N) A^* A (I-\mathcal{R}_N)$ is unitarily equivalent to the semi-infinite section $m,n \geq N$ of the semi-infinite matrix $K$ in equation (\ref{eq:Hankel}).

\bigskip

The largest eigenvalue of this semi-infinite section is further bounded by the Hilbert-Schmidt (Frobenius) norm,
\[
\lambda_1((I-\mathcal{R}_N) A^* A (I-\mathcal{R}_N)) \leq \left[ \sum_{m,n \geq N} K_{m,n}^2 \right]^{1/2}. 
\]
By elementary majorations involving geometric series, there exists $C > 0$ such that the quantity above is less than
\[
C \, \Delta \, N^{-2} \frac{(\rho^*)^{2N}}{1 - (\rho^*)^2}.
\]
For $\rho^* < 1$, it follows that the $N$th singular value of $A_{M,N}$ obeys ($C$ is a number that changes from line to line)
\[
\sigma_N(A_{M,N}) \leq C \, \sqrt{\Delta} \, N^{-1} \, \frac{(\rho^*)^{N}}{\sqrt{1 - (\rho^*)^2}}.
\]

\bigskip

As for the first singular value, we use the assumption that the discretization is reasonable for $A$ (Definition \ref{def:reasonable}) to obtain $\sigma_1(A_{M,N}) \geq C \, \sqrt{\Delta} > 0$ where $C$ is independent of $N$. We assemble inequality (\ref{eq:lower1}) by considering that
\[
\kappa(A_{M,N}) = \frac{\sigma_1(A_{M,N})}{\sigma_N(A_{M,N})},
\]
and noticing that $\rho^* = 1/\alpha$.

\bigskip

So far we have assumed $\alpha < 1 \iff \psi > 0$, but it is clear that the result is also true (and somewhat uninformative) when $\alpha = 1$. The justification of this fact is a very special case of the analysis in the coming section.

\subsection{Fine lower bound (\ref{eq:lower2}) on the condition number}\label{sec:fine-lower}

In this section we consider the worst case scenario when $p^* = \underline{p}$, or equivalently $\rho^* = 1$, or $\psi = 0$. A fortiori the bounds we derive here also hold for any $0 < \rho^* < 1$.

\bigskip

The proof idea for (\ref{eq:lower2}) is that the interval $[\underline{x}, \overline{x}]$ can be subdivided into subintervals of the form $I_j = [\underline{x}, \overline{x}] \cap [(1 + \delta^{j+1})\underline{x}, (1 + \delta^{j}) \underline{x}]$, with $j = J, J+1, \ldots$ and for some $\delta < 1$. Here $J$ is the largest integer such that $(1 + \delta^J) \underline{x} \geq \overline{x}$. The operator $A$ correspondingly splits into the sequence of operators $A_j f = A \chi_{I_j} f$. In accordance with the notation for matrix multiplication we suggestively write $A = (A_J, \ldots, A_j, \ldots )$.

\bigskip

The coarse bound (\ref{eq:lower1}) can now be applied to each $A_j$. The same reasoning as in the previous section applies, yielding
\[
\Sigma_j = \underline{x} \left( 1 + \frac{\delta^j + \delta^{j+1}}{2} \right), \qquad \Delta_j = \underline{x} \left( \frac{\delta^j - \delta^{j+1}}{2} \right),
\]
\[
\tilde{x}_{*,j} = \frac{\delta_j - \underline{x}}{\Delta_j} = \frac{1 + \delta}{1-\delta},
\]
\[
\rho^*_j = \tilde{x}_{*,j} - \sqrt{(\tilde{x}_{*,j})^2 - 1} = \frac{1 - \sqrt{\delta}}{1 + \sqrt{\delta}}.
\]
Hence the sequence of singular values of $A_j$ obeys
\begin{align*}
\sigma_N(A_j) &\leq C \sqrt{\Delta_j} N^{-1} \frac{(\rho^*_j)^N}{\sqrt{1 - (\rho^*_j)^2}}, \\
&\leq D \; (\sqrt{\delta})^j \left( \frac{1 - \sqrt{\delta}}{1 + \sqrt{\delta}} \right)^{N},
\end{align*}
where $D$ is some re-usable constant which depends on $\delta$ and $\underline{x}$, but not $j$ and $N$. For short we let $\eta = \frac{1 - \sqrt{\delta}}{1 + \sqrt{\delta}}$.

\bigskip

The recombination of these various sequences, indexed by $j$, is heuristically done by concatenation. The precise statement is the following inequality due to Weyl.

\begin{lemma} (Weyl)
Consider partitioning a compact operator $A$ as $( B, C )$. The singular values of $A, B$ and $C$ are related by
\[
\sigma^2_{i+j+1}(A) \leq \sigma^2_{i+1}(B) + \sigma^2_{j+1}(C), \qquad i, j \geq 0.
\]
\end{lemma}
\begin{proof}
Write $A A^T = B B^T + C C^T$. Apply Weyl's inequality to this sum of Hermitian compact operators \cite{Weyl, KnutsonTao}:
\[
\lambda_{i+j+1}(A A^T) \leq \lambda_{i+1}(B B^T) + \lambda_{j+1}(C C^T), \qquad i, j \geq 0.
\]
The eigenvalues of $A A^T$, $B B^T$, $C C^T$ are the squares of the singular values of $A$, $B$, $C$ respectively.
\end{proof}

\bigskip

Let us apply this inequality recursively. Let $K > J$ (to be determined) and $n_j$ be integers such that $\sum_{j=J}^K n_j = N$. Then
\[
\sigma_{N+1}^2(A) \leq \sum_{j= J}^K \sigma^2_{n_j+1}(A_j) + \sum_{j=K+1}^\infty \sigma_1^2(A_j).
\]
The last term is seen to be
\[
\sum_{j=K+1}^\infty \sigma_1^2(A_j) \leq D \, \delta^K.
\] 

The numbers $n_j$ are chosen so that each term $\sigma^2_{n_j+1}(A_j) \lesssim \delta^j \eta^{2n_j}$ is also on the order of $\delta^K$ (up to a multiplicative constant that depends on $\delta$, $\underline{x}$, but neither $j$ nor $K$.) For this purpose it is sufficient to take
\[
n_j = \frac{K-j}{2} \log_\eta \delta,
\]
rounded off to the nearest smaller integer.

The sequence $(n_J, \ldots, n_K)$ sums up to a number less than or equal to $N$ provided
\[
N \geq \frac{(K - J + 1)^2}{4} \, \log_\eta \delta.
\]
Choosing $K$ the largest integer smaller than $J - 1 + 2 \sqrt{N / \log_\eta \delta}$ will do. As a result,
\[
\sigma_{N}^2(A) \leq D \,  K  \, \delta^K \leq D \, \sqrt{N} \, \left( \delta^{2/\sqrt{\log_\eta \delta}} \right)^{\sqrt{N}}.
\]
Here and earlier, the proportionality constant $D$ depends on $\delta$ and $\underline{x}$, but not $N$. 

\bigskip

We now address the choice of $0 < \delta < 1$. The number put to the power $\sqrt{N}$ above has for logarithm
\[
\log \left( \delta^{2/\sqrt{\log_\eta \delta}} \right) = 2 \sqrt{\log \eta \, \log \delta},
\]
with $\eta = \frac{1 - \sqrt{\delta}}{1 + \sqrt{\delta}}$. The sharpest bound is obtained when $\log \eta \log \delta$ is minimized as a function of $\delta \in [0,1]$. Numerically, this happens when
\[
\delta = 0.1716...
\]
In that case
\[
\sigma_N(A) \leq D \, N^{1/4} \, \left( \delta^{1/\sqrt{\log_\eta \delta}} \right)^{\sqrt{N}}, \qquad \delta^{1/\sqrt{\log_\eta \delta}} = 0.2875... = e^{-1.2465...}
\]
The largest singular value is lower-bounded away from zero for the same reason as in the previous section. The root-exponential bound on the condition number follows.

\subsection{Upper bound (\ref{eq:upper}) on the condition number}\label{sec:upper}

In this section we seek an upper bound on $\lambda_1(A_{MN}^* A_{MN})$, and a lower bound on $\lambda_N(A_{MN}^* A_{MN})$, for some particular choice of $\mathcal{P}_M$ and $\mathcal{Q}_{N}$. The bound on $\lambda_1$ is easy: use Lemma \ref{lem:one} for $\mathcal{P}_M = I$ and $j=1$ to obtain a bound $C \, \Delta$, independent of $N$. 

\bigskip

$\mathcal{Q}_N$ is chosen as the orthogonal projector in $L^2(\underline{x},\overline{x})$ on the (rescaled Legendre) polynomials of degree $N-1$. For the definition of the orthogonal projector it makes no difference whether those polynomials are orthogonalized or not. Since $x = p^2$, $\mathcal{Q}_N$ is as described in the wording of Theorem \ref{teo:main}. The resulting matrix $\mathcal{Q}_N A^* \mathcal{P}_M A \mathcal{Q}_N$ is analogous to a finite section of $K$, except for the presence of $\mathcal{P}_M$:
\begin{align*}
&\left[ \mathcal{Q}_N A^* \mathcal{P}_M A \mathcal{Q}_N \right]_{m,n} = 2 \, \Delta \, \left( (n+\frac{1}{2})(m+\frac{1}{2}) \right)^{-1/2} \times \\
& \qquad\qquad \int_{\tilde{x}_*}^{\tilde{x}^*} \left[ \, \tilde{x} - \sqrt{\tilde{x}^2 - 1} \,\, \right]^{m+1/2} \; \mathcal{P}_M \left[ \, \tilde{x} - \sqrt{\tilde{x}^2 - 1} \,\, \right]^{n+1/2}  \, d\tilde{x}, \qquad  0 \leq m, n \leq N-1.
\end{align*} 
This expression reduces to $K_{m,n}$ in (\ref{eq:Hankel}) by choosing $M = N$, and $\mathcal{P}_M$ the orthogonal projector in $L^2(x_*,x^*)$, on the subspace
\[
\mbox{span} \{ \left[ \, \tilde{x} - \sqrt{\tilde{x}^2 - 1} \,\, \right]^{n+1/2}  \, \mbox{s.t.} \quad  \tilde{x} = \frac{x - \Sigma}{\Delta}, \, 0 \leq n \leq N-1\}.
\]

\bigskip

We are thus left with the problem of finding a lower bound on the smallest eigenvalue of each finite section $0 \leq m, n \leq N-1$ of the nearly-Hankel matrix $K_{m,n}$ in (\ref{eq:Hankel}). This question was settled in the Hankel case by Szeg\H{o} in 1936 \cite{Szego-1936}, where the full asymptotic behavior as $N \to \infty$ was studied. Widom and Wilf, unaware of Szeg\H{o}'s result, rediscovered it in 1966 with the same techniques \cite{WidomWilf}. Since our matrix $K$ is not exactly of Hankel type (because of the factor $\left( (n+\frac{1}{2})(m+\frac{1}{2}) \right)^{-1/2}$), we rehearse and adapt their argument.

\bigskip

We start with a beautiful characterization of the inverse of a moment matrix which, according to Berg and Szwarc \cite{Berg}, was first discovered by Aitken \cite{Collar}. Let 
\[
H_{m,n} = \int {\rho}^{m+n} d\mu(\rho), \qquad 0 \leq m,n \leq N-1
\]
be a Hankel matrix of moments of the positive measure $\mu(\rho)$. Let $L_n(x)$ denote the orthogonal polynomials associated with $\mu(x)$. Then $H^{-1}$ is similar to the matrix $G$ with entries
\[
G_{m,n} = \frac{1}{2 \pi} \int_0^{2 \pi} L_m(e^{i \theta}) L^*_n(e^{i \theta}) \, d\theta, \qquad 0 \leq m,n \leq N-1.
\]

The main observation of Szeg\H{o}, and Widom and Wilf, is that the large $n$ asymptotics of the polynomials $L_n(z)$ on the unit circle translates into the large $(m,n)$ asymptotics for $G_{m,n}$.

\begin{lemma} (Szeg\H{o}-Widom-Wilf)
Assume that supp $ \mu$ is a finite interval $[\rho_*, \rho^*] \subset \R_{+}$, and that $\mu$ does not vanish on its support\footnote{Szeg\H{o} requires the weaker condition
\[
\int_{\rho_*}^{\rho^*} \frac{\log \mu'(\rho)}{(\rho - \rho^*)^{1/2} (\rho^* - \rho)^{1/2}} , dx < \infty.
\]}. Then
\[
G_{m,n} = \gamma (m+n)^{-1/2} \beta^{m+n} + o((m+n)^{-1/2} \beta^{m+n}),
\]
for some $\gamma > 0$, and where
\begin{equation}\label{eq:betarho}
\beta = \frac{\rho^* + \rho_* + 2}{\rho^* - \rho_*} + \left[ \left( \frac{\rho^* + \rho_* + 2}{\rho^* - \rho_*} \right)^2 - 1 \right]^{1/2}.
\end{equation}
\end{lemma}

The next step is to approximate the eigenvector corresponding to the leading eigenvalue of $G$. This is where we depart from \cite{Szego-1936, WidomWilf}. Our matrix of interest is not $H^{-1}$ but
\[
(K^{-1})_{m,n} \; \mbox{similar to} \; \frac{1}{\Delta} (m+\frac{1}{2})^{1/2} G_{m,n} (n+\frac{1}{2})^{1/2}.
\]
Consider approximating the finite section $0 \leq m,n \leq N-1$ as
\[
\frac{1}{\Delta} (m+\frac{1}{2})^{1/2} G_{m,n} (n+\frac{1}{2})^{1/2} = L_{m,n} + R_{m,n},
\]
where $L_{m,n}$ is the leading rank-1 expression
\[
L_{m,n} = \frac{\gamma}{2 \Delta} (2N-2)^{-1/2} v_m v_n, \qquad v_n = (n+\frac{1}{2})^{1/2} \beta^n.
\]
It is easy to show that the spectral radius of the remainder $R_{m,n}$ tends to zero as $N \to \infty$, so it can be neglected in an asymptotic sense for large $N$. On the other hand, $v_n$ is the eigenvector that corresponds to the unique nonzero eigenvalue of the leading part $L_{m,n}$. This eigenvalue obeys
\[
\lambda_1(L) = \frac{\gamma}{\Delta} (2N-2)^{-1/2} \sum_{n=0}^{N-1} |v_n|^2 \leq C \frac{\sqrt{N}}{\Delta} \frac{\beta^{2N}}{\beta^2-1}. 
\]
By taking the constant sufficiently large, this bound also holds for $L + R$, for all $N$.  Specializing to $K_{m,n}$, we get
\[
\lambda_N(K) \geq C \frac{\Delta}{\sqrt{N}} (\beta^2 - 1) \beta^{-2N},
\]
with $\beta$ given in (\ref{eq:betarho}). The result on the condition number follows.

\subsection{Szeg\H{o} average decay of the eigenvalues.}\label{sec:average}

New ideas may be required to sharpen the constants in the lower bounds on the condition number. One useful piece of information could be the rate at which the determinant of a Hankel form grows as $N \to \infty$.

\bigskip

Consider $D_N = \det H^{(N+1)}$, where
\[
H^{(N+1)}_{m,n} = \int_{\rho_*}^{\rho^*} \rho^{m+n} d\mu(\rho), \qquad 0 \leq m,n \leq N.
\]
Szeg\H{o}\footnote{Explained on p.85 of \cite{Szego-book}. Szeg\H{o} is best known for proving the corresponding result for Toeplitz forms (in which case $c=1$) when he was an undergraduate student, after P\'{o}lya posed it as a conjecture. This note on the historical context of the Szeg\H{o} distribution theorem is taken from \cite{Simon}.} found the asymptotic expression
\[
\lim_{n \to \infty} c^{-2N-1} \frac{D_N}{D_{N-1}} = 2 \pi \exp \left( \frac{1}{2 \pi} \int_0^{2 \pi} \log \mu'(h(\theta)) d\theta \right),
\]
where $h(\theta)$ is a map from $[0,2 \pi)$ to $[\rho_*, \rho^*]$, and $c$ is called the transfinite diameter of $[\rho_*, \rho^*]$ corresponding to this map. It is possible to choose $h$ as a simple trigonometric function such that
\[
c = \frac{\rho^* - \rho_*}{4}.
\]
If we let $R$ for the right-hand side, we obtain the explicit asymptotic formula
\begin{align*}
D_N &\sim R^N \left[ c^{2N+1} \,  c^{2N-1} \ldots c^3 \, c \right] D_0, \\
&\sim R^N c^{(N+1)^2} D_0.
\end{align*}
Since $D_N$ is the product of the $N+1$ eigenvalues of $H^{(N+1)}$, and \emph{if} we postulate that these eigenvalues decay geometrically, then the only possible decay rate is (up to a polynomial factor) $c^{n}$ with $c = \frac{\rho^* - \rho_*}{4}$. The same decay rate would hold for the eigenvalues of the matrix $K$ in (\ref{eq:Hankel}). The corresponding asymptotics for $\kappa_N$ for the finite-section discretization of the operator $A$ would follow as $\kappa_N \sim c^{-N/2}$ up to a polynomial factor.

\bigskip

As we have no indication that the eigenvalues of sections of $H$ or $K$ indeed decay geometrically, or whether the ``average" rate $c^{-N/2}$ could be useful in any way toward formulating a bound on $\kappa_N$, we contend ourselves with reporting it numerically with the other bounds in Section \ref{sec:num}.

\section{Small $p$ asymptotics}\label{sec:asymp}

The obstruction to the traveltime tomography problem in the case $p = 0$ (rays perpendicular to the layering) was covered in Section \ref{sec:asymp-intro}. In the more general case when $0 \leq p \leq p^*$ with small $p^*$, equation (\ref{eq:tau}) can be rewritten at $z = h$ as
\[
\tau(h,p) = \int_0^h \frac{1}{c(z) \sqrt{1 - p^2 c^2(z)}} \, dz.
\]
Perform a binomial expansion of the inverse square root to find its Taylor expansion as
\[
(1 + x)^{-1/2} = \sum_{n = 0}^\infty \begin{pmatrix} -1/2 \\ n \end{pmatrix} x^n, \qquad \mbox{if } |x| < 1.
\]
The first few generalized binomial coefficients are
\[
\begin{pmatrix} -1/2 \\ 0 \end{pmatrix} = 1, \qquad \begin{pmatrix} -1/2 \\ 1 \end{pmatrix} = - \frac{1}{2}, \qquad \begin{pmatrix} -1/2 \\ 2 \end{pmatrix} = \frac{3}{8}, \; \mbox{etc.}
\]
In our case, if $p$ is small enough $|x| = p^2 c^2(z) < 1$. The smaller $p$ the more accurate the truncation of the sum to the first few terms:
\[
\tau(h,p) = \sum_{n=0}^{d-1} \begin{pmatrix} -1/2 \\ n \end{pmatrix} (-1)^n p^{2n} \int_0^h c^{2n-1}(z) \, dz+ O((c_0 p)^{2d}).
\]
(We placed the ad-hoc factor $c_0 = c(z=0)$ in the remainder to make it dimensionless. Recall that $c_0 p = \cos \theta$ where $\theta$ is the angle that the ray labeled $p$ makes with the surface $z=0$.) It is now clear that if two profiles $c_1$ and $c_2$ have matching ``odd moments" up to degree $d-1$, i.e.
\[
\int_0^h c_1^{2n-1}(z) \, dz = \int_0^h c_2^{2n-1}(z) \, dz, \qquad 0 \leq n < d,
\]
their responses $\tau_1(h,p)$ and $\tau_2(h,p)$ will match to within $O((c_0 p)^{2d})$. The name ``moment" owes from the fact that these integrals are precisely moments of the slowness distribution function introduced in equation (\ref{eq:SDF}), namely
\[
\int_0^h c^{m}(z) \, dz = \int q^{m} F(q;p) \ dq.
\]
The moment-matching inverse problem is notoriously ill-posed \cite{ShohatTamarkin}.

\bigskip

The reasoning carries over without difficulty to the case of slightly differing odd moments. For instance, in order to get
\[
|\tau_1(h,p) - \tau_2(h,p)| \lesssim \eps, \qquad \mbox{for } 0 \leq p \leq p^*,
\]
it suffices to find the smallest $d$ such that ${(c_0 p^*)}^{2d} \leq \eps$, and require
\[
\Bigl\lvert \int_0^h c_1^{2n-1}(z) \, dz - \int_0^h c_2^{2n-1}(z) \, dz \Bigr\rvert \leq \eps \, \frac{p^{-2n}}{d \begin{pmatrix} -1/2 \\ n \end{pmatrix}}, \qquad \mbox{for } 0 \leq n < d.
\]

\bigskip

This latter relation defines a rather elongated set around $c_1$, a ``ambiguity region" of kinematically near-equivalent velocity profiles $c_2$. In general there will exist such near-equivalent $c_1$ and $c_2$ for which the difference $c_2 - c_1$ is non-oscillatory, i.e., contain only low wavenumbers. A numerical illustration of this phenomenon is shown in Section \ref{sec:num}. As a result \emph{Tychonov regularization will hardly be able to discriminate between $c_1$ and $c_2$} if they are comparably smooth: this is bad news for the prospect of solving the inverse problem.

\section{Numerics}\label{sec:num}

In Figure \ref{fig:egval-abel} we show an illustration of the various bounds on the condition number $\kappa_N$ as a function of $N$. For the particular choice of discretization made for the upper bound (\ref{eq:upper}) in Theorem \ref{teo:main}, recall that 
\[
\kappa_N = \frac{\sigma_1}{\sigma_N} = \sqrt{\frac{\lambda^{(N)}_{1}}{\lambda^{(N)}_{N}}},
\]
where $\lambda^{(N)}_{n}$ is the $n$th eigenvalues of the size-$N$ finite section $\{ K_{m,n} : 0 \leq m,n < N \}$ of the infinite matrix $K$ in (\ref{eq:Hankel}). The graphs of $1/\sqrt{\lambda^{(N)}_{n}}$ are plotted on the same picture as a function of $n$, with the different curves indexed by $N$. Other discretization choices may not be linked in any way to the $\lambda^{(N)}_{n}$.

\bigskip

\subsection{Bounds on the condition number}
\begin{figure}[H]
\centering
\includegraphics[width=8cm]{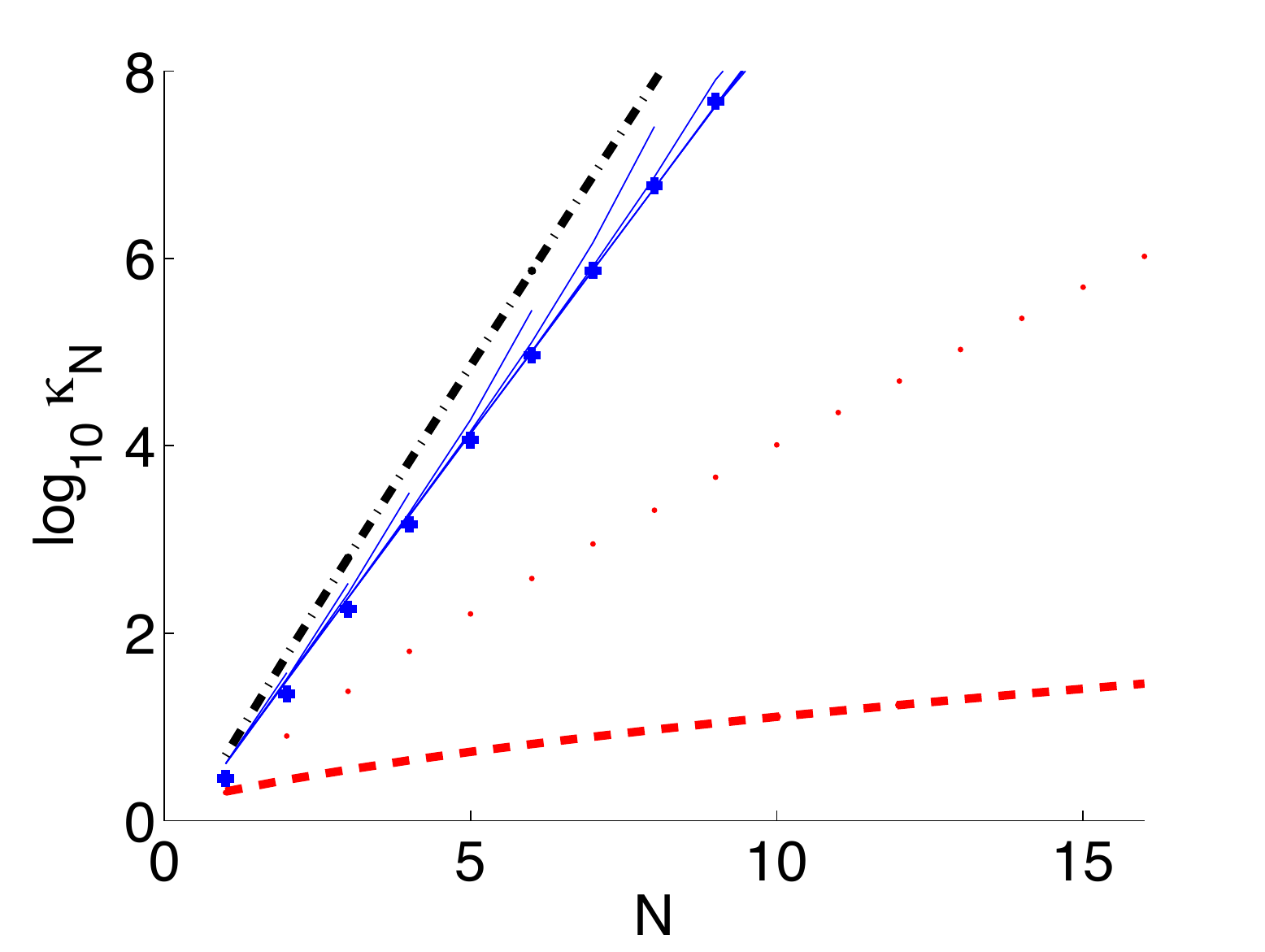}
\hfill
\includegraphics[width=8cm]{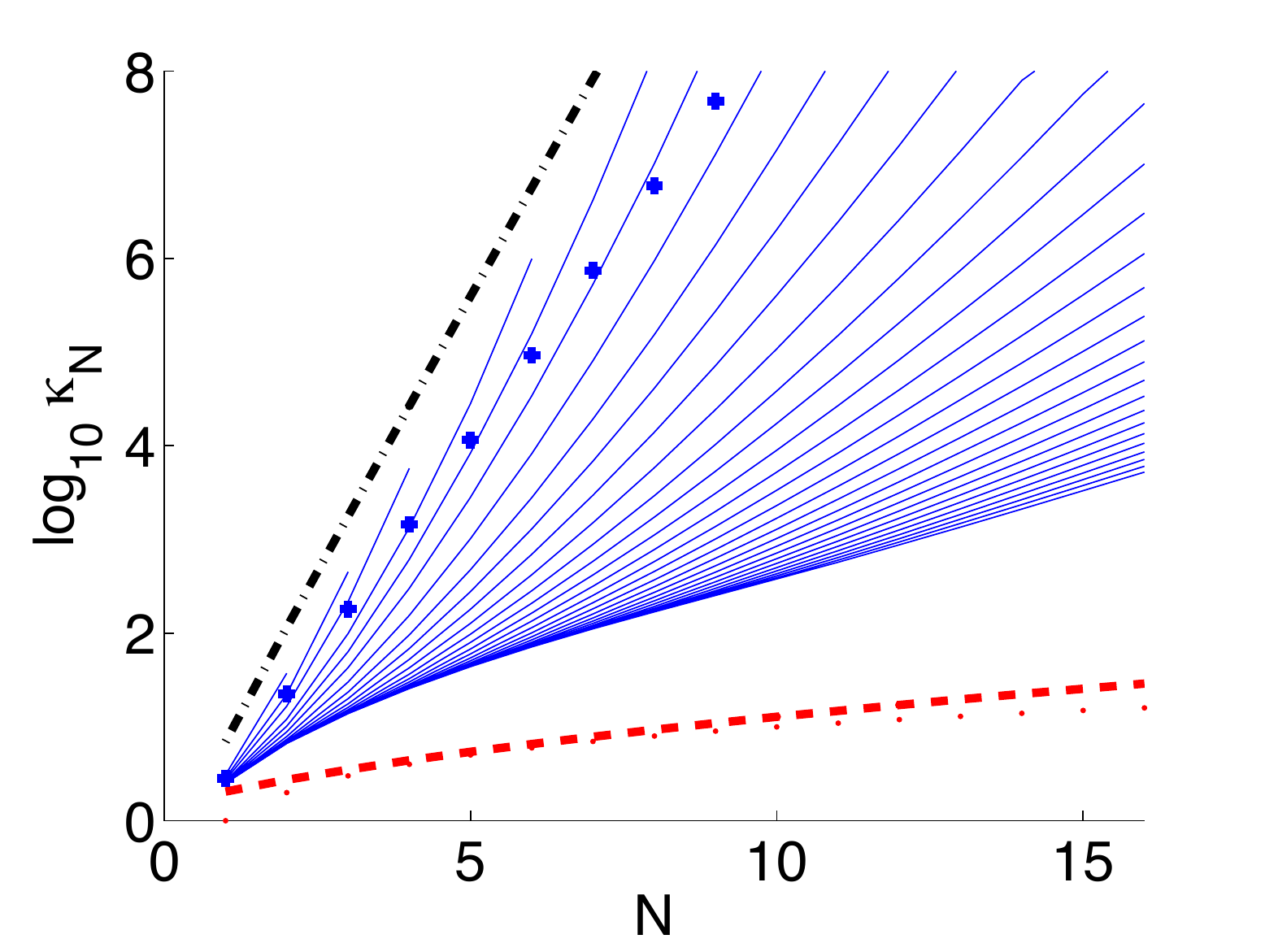}
\caption{Various bounds on the condition number as a function of the discretization parameter $N$, for the reflected rays setup. Notice the logarithmic scale of the $y$ axis. Left: $\rho_* = 0$ and $\rho^* = 0.5$. Right: $\rho_* = 0.5$ and $\rho^* = 1$. Dotted curve: first lower bound (\ref{eq:lower1}). Dashed curve: second lower bound (\ref{eq:lower2}). Dash-dotted curve: upper bound (\ref{eq:upper}). Solid curves: plots of $1/\sqrt{\lambda^{(N)}_{n}}$ (see text), as a function of $1 \leq n \leq N$ for varying $N$. Blue crosses: Szeg\H{o} average rate from Section \ref{sec:average}. See further comments in the text.}\label{fig:egval-abel}
\end{figure}

\bigskip

The curves for the bounds were scaled by an arbitrary constant, which amounts to an arbitrary vertical translation in logarithmic scale. The first observation is that the upper bound (\ref{eq:upper}) seems sharp as it scales like $1/\sqrt{\lambda^{(N)}_{N}}$. 

\bigskip

The behavior of the eigenvalues $\lambda^{(\infty)}_{n}$ of the infinite matrix $K$ is given by the lower envelope of the eigenvalue curves. Its scaling seems to be root-exponential in the case $\rho^* = 1$ (Figure \ref{fig:egval-abel}, right panel), i.e., of the form $c e^{-d\sqrt{n}}$ for some numbers $c, d > 0$.  The lower bound (\ref{eq:lower2}) indeed scales root-exponentially, albeit with a different non-sharp constant in the exponential. Note that the ratio $\sqrt{\lambda^{(\infty)}_{1} / \lambda^{(\infty)}_{n}}$ in the case $N \to \infty$ can be seen as the \emph{discretization-free condition number}, i.e., the condition number of the best discretization which gives rise to the largest singular values. In that case, the projector $\mathcal{Q}_n$ project onto the subspaces formed by the $n$ eigenvectors corresponding to the largest eigenvalues $\lambda^{(\infty)}_{1}, \ldots, \lambda^{(\infty)}_{n}$. The discrepancy between $\lambda^{(N)}_{N}$ (upper end of the curves) and $\lambda^{(\infty)}_{N}$ (lower envelope at the same abscissa) shows that the discretization defined by taking finite sections is far from being ``best" in the sense discussed above.

\subsection{Negative implications for imaging}

\begin{figure}[H]
\begin{center}
\leavevmode
\includegraphics[width=\textwidth]{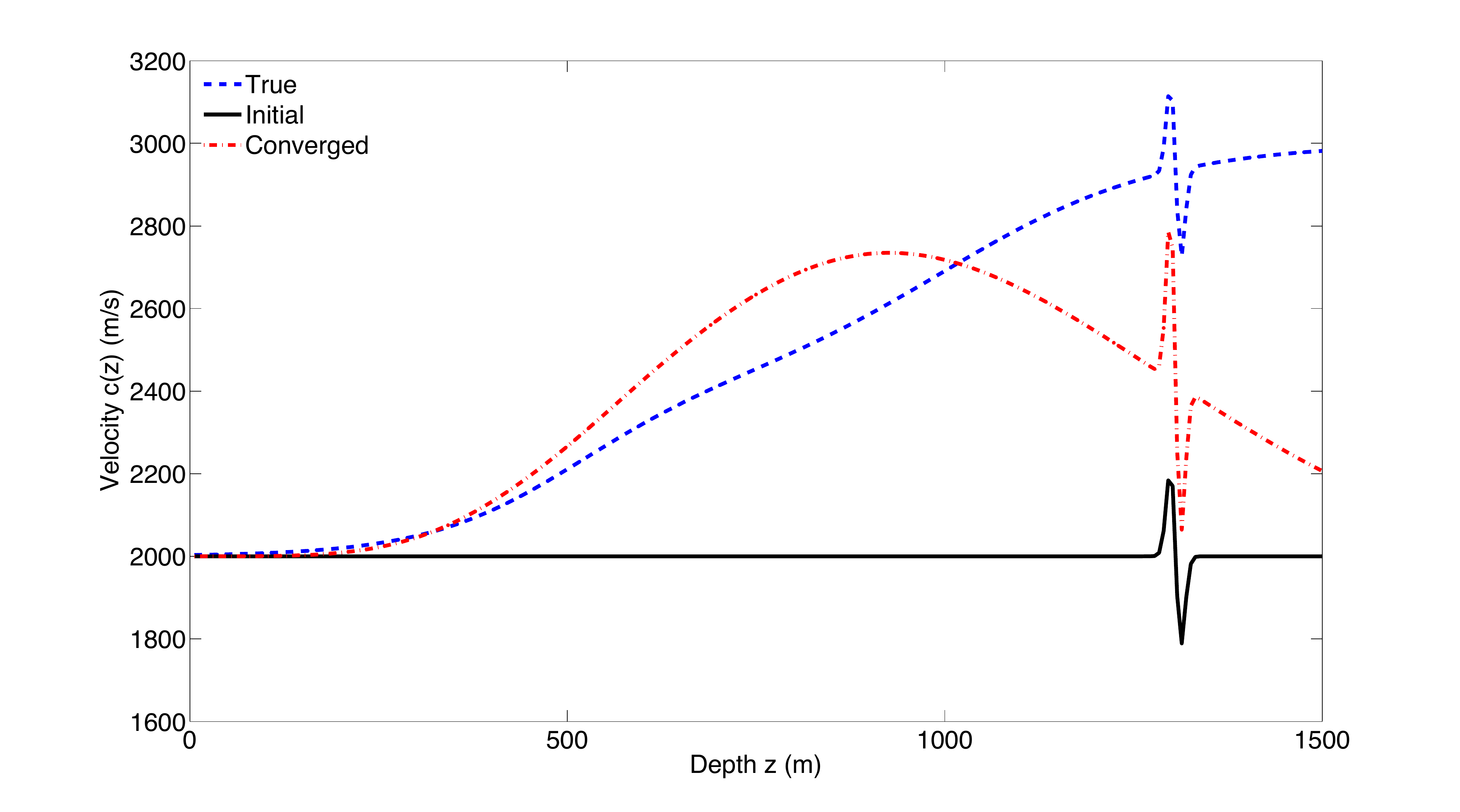}
\end{center}
\caption{
True vs. converged velocities. Velocity profiles are functions of $z$ only. Different levels of Tychonov regularization do not noticeably improve the converged model.}
\label{fig:c(z)}
\end{figure}

Although the theory in this paper concerns rays rather than waves, the conditioning issue identified here also plagues the finite-frequency, waveform-based inversion problem of reflection seismology. A finite difference acoustic wave simulation was carried out in a smoothly increasing medium $c(z)$ shown as the blue dashed curve in Figure \ref{fig:c(z)}, to create synthetic seismograms of reflected waves (not shown). The receivers cover the surface $z=0$. There is a single source at $x = z = 0$. Note the ``reflector" near $z \simeq 1300$ which is responsible for the wave echos recorded at the surface. The wavelength of the probing waves is about 50 m. The data corresponding to diving waves are discarded. The initial $c(z)$ is the black solid curve.

\bigskip

The inverse problem of determining the background velocity $c(z)$ from these synthetic seismograms was solved using (our own implementation of) the Mulder-Van Leeuwen correlation-focusing method \cite{Mulder}. In a nutshell, least-squares based inversion -- minimizing the $\ell_2$ norm of the waveform residual -- would fail because of lack of convexity of the minimization objective, but correlation focusing is an alternative choice of objective that does not suffer (as much) from that problem. 

\bigskip

\emph{The inversion procedure converges successfully so that data are fit within a few digits of accuracy. Yet the converged speed profile (red dash-dotted curve) is significantly different from the original ``true" speed profile. Various levels of Tychonov regularization did not help in improving convergence.}

\bigskip

\bigskip

The original $c(z)$ used in the forward modeling step was chosen to increase monotonically, so the ill-conditioning is not due to the rearrangement ambiguity (also called ``presence of low-velocity zones".)  Instead, it is a (spectacular) finite-frequency remnant of the conditioning problem associated with reflected rays as studied in this paper.

\bigskip

\begin{figure}[H]
\begin{center}
\includegraphics[width=17cm]{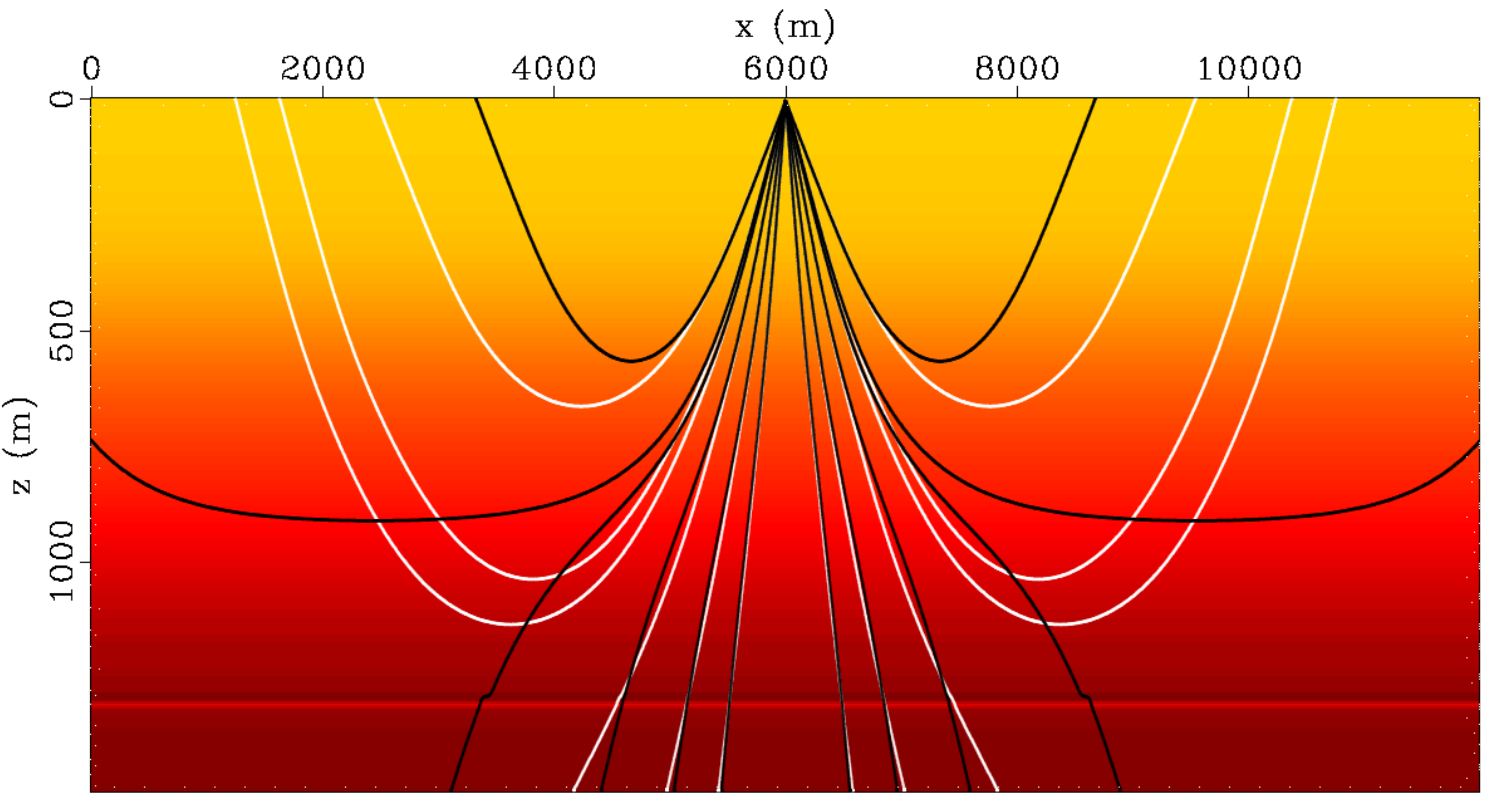}
\end{center}
\caption{
White and black lines are rays in the true and converged background velocity shown in Figure \ref{fig:c(z)}, respectively.
The number of lines and take-off angles are same in both cases.
Yellow and red colors represent the magnitude of the true background velocity model. The authors of RSF and Madagascar are gratefully acknowledged for providing the plotting routines.
}
\label{fig:rays-tomog}
\end{figure}

Another numerical piece of evidence for the problem associated with reflected rays is Figure \ref{fig:rays-tomog}. The background velocity is shown in shades of yellow and red: it is the same ``true" wave speed profile as earlier. In white, the rays of geometrical optics were traced in this ``true" medium (blue dashed curve in Figure \ref{fig:c(z)}). In black, we traced rays in the converged ``optimal" medium from correlation-focusing inversion (red dash-dotted curve in Figure \ref{fig:c(z)}.) Notice how the transmitted rays reach the reflector $z \simeq 1300$ (and then reflect) at almost the same location as the white rays with the same take-off angle. The diving rays are completely different, on the other hand -- hence they contain much more information than the reflected rays.

\bigskip

Finally, we compare the odd moments of the converged velocity profile $c_k(z)$ from the correlation-focusing method after $k$ iterations, to those of the ``true" velocity profile $c(z)$. The table below lists the quantity
\[
\frac{\int (c(z))^p - (c_k(z))^p dz }{\int (c(z))^p dz }
\]
for different values of $p$ and $k$.

\begin{table}[H]
\begin{center}
\label{table:moment}
\begin{tabular}{l|llll}
$p$ $\backslash$ $k$ &   0 &      1    &      2    &     3     \\
\hline
-1  &  -1.66e-01 & 2.51e-02&  6.20e-03 &  6.62e-03   \\
1   & 1.56e-01 & -2.72e-02 & -4.84e-03 & -5.32e-03    \\
3   & 4.27e-01 & -8.43e-02 & -7.41e-03 & -9.00e-03    \\
5   & 6.34e-01 & -1.38e-01 &  3.25e-03 &  4.30e-04    \\
7   & 7.78e-01 & -1.82e-01 &  2.88e-02 &  2.47e-02    \\
9   & 8.71e-01 & -2.12e-01 &  6.80e-02 &  6.28e-02  
\end{tabular}
\end{center}
\end{table}

The moments match to within a few digits after very few iterations, as they should from the discussion in Section \ref{sec:asymp}.

\section{Discussion}

We have shown that the isotropic, laterally-homogeneous traveltime tomography inverse problem has well-posed formulations in the case of diving rays, but suffers from incurable ill-conditioning in the case of reflected rays. While diving rays involve a Volterra integral equation, reflected rays involve a Fredholm integral equation.  Intuitively, a Fredholm operator is to a rank-deficient full matrix what a Volterra operator is to the upper-triangular restriction of such a matrix.

\bigskip

Our analysis shows that well-posedness is linked to the presence of overturning rays, i.e., rays whose direction is at some point parallel to the level lines of the speed profile. We do not know if this non-transversality condition could play a role for the analysis of the more general case of a laterally varying $c(x,z)$. 

\bigskip

The ill-conditioned nature of the reflection traveltime tomography problem has serious implications for imaging, even at finite frequencies. The seismic inverse problem in a smooth, layered background $c(z)$ with surface data can only be only well-posed if either (1) low-frequency data is seriously taken into account, and/or (2) the reflectors are more or less ``dense" in the sense that the true wave speed profile is ``rough everywhere". The latter point was made precise by Symes who wrote a remarkable stability estimate in \cite{Symes}.

\bigskip

Finally, it should be mentioned that if we restrict the domain to a rectangle, and avail ourselves of complete data on all the sides, then the problem of recovering the wave speed from traveltime data becomes much better posed. For instance, Mukhometov proved a stability estimate (with loss of one derivative) in the case of isotropic media that deviate little from a constant \cite{Mukhometov}.

\end{document}